%%%%%%%%%%%%%%%%%%%%%%%%%%%%%%%%
%
% These are the macroes
%
%%%%%%%%%%%%%%%%%%%%%%%%%%%%%%%%

\magnification\magstephalf

\voffset0truecm
\hoffset=0truecm
\vsize=23truecm
\hsize=15.8truecm
\topskip=1truecm

\binoppenalty=10000
\relpenalty=10000

%%%% Famiglie di Fonts
\font\tenbb=msbm10		\font\sevenbb=msbm7		\font\fivebb=msbm5
\font\tensc=cmcsc10		\font\sevensc=cmcsc7 	\font\fivesc=cmcsc5
\font\tensf=cmss10		\font\sevensf=cmss7		\font\fivesf=cmss5
\font\tenfr=eufm10		\font\sevenfr=eufm7		\font\fivefr=eufm5

%\font\eurb=eurb10

\newfam\bbfam	\newfam\scfam	\newfam\frfam	\newfam\sffam

\textfont\bbfam=\tenbb
\scriptfont\bbfam=\sevenbb
\scriptscriptfont\bbfam=\fivebb

\textfont\scfam=\tensc
\scriptfont\scfam=\sevensc
\scriptscriptfont\scfam=\fivesc

\textfont\frfam=\tenfr
\scriptfont\frfam=\sevenfr
\scriptscriptfont\frfam=\fivefr

\textfont\sffam=\tensf
\scriptfont\sffam=\sevensf
\scriptscriptfont\sffam=\fivesf

\def\bb{\fam\bbfam \tenbb} % Blackboard bold
\def\sc{\fam\scfam \tensc} % Maiuscoletto
 % Gotico
 % Sans serif

%%%% Fonts

\font\sezfont=cmbx10 scaled \magstep1
\font\subsectfont=cmbx10 scaled \magstephalf
\font\titfont=cmbx10 scaled \magstep2
\font\autfont=cmcsc10
\font\intfont=cmss10 

%%%% Abbreviazioni
\let\no=\noindent
\let\bi=\bigskip
\let\me=\medskip
\let\sm=\smallskip
\let\ce=\centerline

\let\io=\infty
\def\qqquad{\quad\qquad}

%%%% Numerazione

\newcount\sectno\sectno=0
\newcount\subsectno\subsectno=0
\newcount\thmno\thmno=0
\newcount\tagno\tagno=0
\newcount\notitolo\notitolo=0
\newcount\defno\defno=0

% sezione
\def\sect#1\par{
	\global\advance\sectno by 1 \global\subsectno=0\global\defno=0\global\thmno=0
	\vbox{\vskip.75truecm\advance\hsize by 1mm
	\hbox{\centerline{\sezfont \the\sectno.~~#1}}
	\vskip.25truecm}\nobreak}

% sottosezione
\def\subsect#1\par{
	\global\advance\subsectno by 1
	\vbox{\vskip.75truecm\advance\hsize by 1mm
	\line{\subsectfont \the\sectno.\the\subsectno~~#1\hfill}
	\vskip.25truecm}\nobreak}
	
% definizioni
\def\defin#1{\global\advance\defno by 1
	\global\expandafter\edef\csname+#1\endcsname%
    {\number\sectno.\number\defno}
    \no{\bf Definition~\the\sectno.\the\defno.}}

% teoremi, lemmi, proposizioni e osservazioni
\def\thm#1#2{
	\global\advance\thmno by 1
	\global\expandafter\edef\csname+#1\endcsname%
	{\number\sectno.\number\thmno}
	\no{\bf #2~\the\sectno.\the\thmno.}}

% equazioni
\def\Tag#1{\global\advance\tagno by 1 {(\the\tagno)}
    \global\expandafter\edef\csname+#1\endcsname%
    		{(\number\tagno)}}
\def\tag#1{\eqno\Tag{#1}}

\def\rf#1{\csname+#1\endcsname\relax}

% fine dimostrazione
\def\proof{\no{\sl Proof.}\enskip}
\def\qedn{\thinspace\null\nobreak\hfill\hbox{\vbox{\kern-.2pt\hrule height.2pt
        depth.2pt\kern-.2pt\kern-.2pt \hbox to2.5mm{\kern-.2pt\vrule
        width.4pt \kern-.2pt\raise2.5mm\vbox to.2pt{}\lower0pt\vtop
        to.2pt{}\hfil\kern-.2pt \vrule
        width.4pt\kern-.2pt}\kern-.2pt\kern-.2pt\hrule height.2pt
        depth.2pt \kern-.2pt}}\par\medbreak}
    \def\qed{\hfill\qedn}
    
%%% Numerazione di pagina / intestazioni
\newif\ifpage\pagefalse
\newif\ifcen\centrue

\headline={
\ifcen\hfil\else
\ifodd\pageno
\global\hoffset=0.5truecm
\else
\global\hoffset=-0.4truecm
\fi\hfil
\fi}

\footline={
	\ifpage
		\hfill\rm\folio\hfill
	\else
		\global\pagetrue\hfill
\fi}

\lccode`\'=`\'

%%% Riferimenti bibliografici
\def\bib#1{\me\item{[#1]\enskip}}

%%% Macrodefinizioni matematiche
\def\ca#1{{\cal #1}}
\def\C{{\bb C}} \def\d{{\rm d}}
\def\R{{\bb R}} \def\Q{{\bb Q}}
\def\Z{{\bb Z}} \def\N{{\bb N}}

\let\eps=\varepsilon \let\phe=\varphi

\mathchardef\void="083F
 % PARTE REALE
 % PARTE IMMAGINARIA

\def\diag{\mathop{\rm Diag}\nolimits}

\def\id{\mathop{\rm Id}\nolimits}

\def\res{{\rm Res}}

\def\invlim{\mathop{\vtop{\offinterlineskip
\hbox{\rm lim}\kern1pt\hbox{\kern-1.5pt$\longleftarrow$}\kern-3pt}
}\limits}
\def\neweq#1$${\xdef #1{(\the\capno.\the\tagno)}
	\eqno #1$$
	\iffinal\else\rsimb#1\fi
	\global \advance \tagno by 1}
\def\neweqa#1$${\xdef #1{{\rm(\the\capno C.\the\tagno)}}
	\eqno #1$$
	\iffinal\else\rsimb#1\fi
	\global \advance \tagno by 1}
\def\newforclose#1{%
	\global\advance\tagno by 1 
    \global\expandafter\edef\csname+#1\endcsname%
    		{(\number\capno.\number\tagno)}
		\hfil\llap{$(\the\capno.\the\tagno)$}\hfilneg}
\def\newforclosea#1{
	\xdef #1{{\rm(\the\capno C.\the\tagno)}}
	\hfil\llap{$#1$}\hfilneg
	\global \advance \tagno by 1
	\iffinal\else\rsimb#1\fi}
\def\forevery#1#2$${\displaylines{\let\neweqa=\newforclosea
	\let\tag=\newforclose\hfilneg\rlap{$\qqquad\forall#1$}\hfil#2\cr}$$}
\def\casi#1{\vcenter{\normalbaselines\mathsurround=0pt
		\ialign{$##\hfil$&\quad##\hfil\crcr#1\crcr}}}

%%%%%%%%%%%%%%%%%%%%%%%%%%%%%%%%%%%%%%%%%%%%%%%%%%%%%%%%%%%%%%%%%%%%%%%%%%%%%%%%%%%%%%%%%%%%%%%%%%%%%%%%%%%%%%%%%%%%%%%%%%%%%%%%%%%%%%%%%%%%%%%%%%%%%%%%%%%%%%%%%%%%%%%%%%%%%%%%%%%%%%%%%%%%%%%%%%%%%%%%%%%%%%%%%%%%%%%%%%%%%%%%%%%%%%%%%%%%%%%%

%%%%%%%%%%%%%%%%%%%%%%%%%%%%%%%%
%
% Begin of the paper
%
%%%%%%%%%%%%%%%%%%%%%%%%%%%%%%%%

\input graphicx.tex
\input xy
\xyoption{all}
%\input macro.tex

%\hfill{\newcount\versiono\versiono=3
%{\tt Version \the\versiono; \the\day/\the\month/\the\year}}

\ce{\titfont Brjuno conditions for linearization } 
\ce{\titfont in presence of resonances} 

%\ce{\titfont }

\me\ce{\autfont Jasmin Raissy}
\sm\ce{\intfont Dipartimento di Matematica, Universit\`a di Pisa}

\ce{\intfont Largo Bruno Pontecorvo 5, 56127 Pisa}

\sm\ce{\intfont E-mail: {\tt raissy@mail.dm.unipi.it}}
\bi

{\narrower

{\sc Abstract.} We present a new proof, under a slightly different (and more natural) arithmetic hypothesis, and using direct computations via power series expansions, of a holomorphic linearization result in presence of resonances originally proved by R\"ussmann.
 
}

\bi

\footnote{}{{\sl \hskip-20pt\noindent Mathematics Subject Classification (2010).} Primary 37F50; Secondary 32H50. \hfill \break {\sl Key words and phrases.} Linearization problem, resonances, Brjuno condition, discrete local holomorphic dynamical systems.}

%%%%%%%%%%%%%%%%%%%%%%%%%%%%%%%%%%%%%%%%%%%%%%%%%%%%%%%%%%%%%%%%%%%%%%%%%%%%%%%%
\sect Introduction
%%%%%%%%%%%%%%%%%%%%%%%%%%%%%%%%%%%%%%%%%%%%%%%%%%%%%%%%%%%%%%%%%%%%%%%%%%%%%%%%

We consider a germ of biholomorphism~$f$ of~$\C^n$ at a fixed point~$p$, which, up to translation, we may place at the origin~$O$. One of the main questions in the study of local holomorphic dynamics (see [A1], [A2], [Bra], or [R3] Chapter 1, for general surveys on this topic) is when~$f$ is {\it holomorphically linearizable}, i.e., when there exists a local holomorphic change of coordinates such that~$f$ is conjugated to its linear part $\Lambda$. 

A way to solve such a problem is to first look for a formal transformation $\phe$ solving 
$$
f\circ\phe = \phe \circ \Lambda,
$$
i.e., to ask when $f$ is {\it formally linearizable}, and then to check whether $\phe$ is convergent. Moreover, since up to linear changes of the coordinates we can always assume $\Lambda$ to be in Jordan normal form, i.e.,
$$
\Lambda = \pmatrix { \lambda_1 & & \cr
		 			\eps_2   & \lambda_2 & \cr
					 &  \ddots & \ddots & \cr
					 &  & \eps_n & \lambda_n  },
$$
where the eigenvalues $\lambda_1, \dots, \lambda_n\in\C^*$ are not necessarily distincts, and $\eps_j\in\{0,\eps\}$ can be non-zero only if $\lambda_{j-1} = \lambda_j$, we can reduce ourselves to study such germs, and to search for $\phe$ {\it tangent to the identity}, that is, with linear part equal to the identity.

\sm The answer to this question depends on the set of eigenvalues of~$\d f_O$, usually called the {\it spectrum} of~$\d f_O$. In fact, if we denote by~$\lambda_1, \dots, \lambda_n\in \C^*$ the eigenvalues of~$\d f_O$, then it may happen that there exists a multi-index~$Q=(q_1, \dots, q_n)\in \N^n$, with~$|Q|\ge 2$, such that
$$
\lambda^Q - \lambda_j:=\lambda_1^{q_1}\cdots\lambda_n^{q_n} - \lambda_j = 0\tag{eqres}
$$
for some~$1\le j\le n$; a relation of this kind is called a {\it (multiplicative) resonance} of~$f$ {\it relative to the $j$-th coordinate}, $Q$ is called a {\it resonant multi-index relative to the $j$-th coordinate}, and we put
$$
{\res}_j(\lambda):=\{Q\in\N^n\mid |Q|\ge 2, \lambda^Q = \lambda_j\}.
$$ 
The elements of $\res(\lambda) := \bigcup_{j=1}^n \res_j(\lambda)$ are simply called {\it resonant}  multi-indices. A {\it resonant monomial} is a monomial~$z^Q:=z_1^{q_1}\cdots z_n^{q_n}$ in the~$j$-th coordinate with $Q\in\res_j(\lambda)$.

\sm Resonances are the formal obstruction to linearization. Indeed, we have the following classical result:

\sm\thm{TePD}{Theorem} (Poincar\'e, 1893 [P]; Dulac, 1904 [D]) {\sl Let~$f$ be a germ of biholomorphism of~$\C^n$ fixing the origin~$O$ with linear part in Jordan normal form. Then there exists a formal transformation $\phe$ of $\C^n$, without constant term and tangent to the identity, conjugating~$f$ to a formal power series~$g\in\C[\![z_1, \dots, z_n]\!]^n$ without constant term, with same linear part and containing only resonant monomials. Moreover, the resonant part of the formal change of coordinates $\psi$ can be chosen arbitrarily, but once this is done, $\psi$ and $g$ are uniquely determined. In particular, if the spectrum of $\d f_O$ has no resonances, $f$ is formally linearizable and the formal linearization is unique.}

\sm A formal transformation $g$ of $\C^n$, without constant term, and with linear part in Jordan normal form with eigenvalues $\lambda_1,\dots,\lambda_n\in\C^*$, is called {\it in Poincar\'e-Dulac normal form} if it contains only resonant monomials with respect to $\lambda_1,\dots,\lambda_n$.  

If $f$ is a germ of biholomorphism of $\C^n$ fixing the origin, a series $g$ in Poincar\'e-Dulac normal form formally conjugated to $f$ is called a {\it Poincar\'e-Dulac (formal) normal form of~$f$}. 

%The formal series~$g$ is called a {\it Poincar\'e-Dulac normal form} of~$f$. 

% From the formal point of view, we have the following classical result (see [Ar] pp.~192--193 for a proof):

%\sm\thm{Te0.1}{Theorem} {\sl Let~$f$ be a germ of holomorphic diffeomorphism of~$\C^n$ fixing the origin~$O$ with no resonances. Then~$f$ is formally conjugated to its differential~$\d f_O$.}

%\sm In presence of resonances, even the formal classification is not easy, as the following result of Poincar\'e-Dulac, [Po], [D], shows

%\sm\thm{Te0.2}{Theorem} (Poincar\'e-Dulac) {\sl Let~$f$ be a germ of holomorphic diffeomorphism of~$\C^n$ fixing the origin~$O$. Then~$f$ is formally conjugated to a formal power series~$g\in\C[\![z_1, \dots, z_n]\!]^n$ without constant term such that~$\d g_O$ is in Jordan normal form, and~$g$ has only resonant monomials.}

%\sm The formal series~$g$ is called a {\it Poincar\'e-Dulac normal form} of~$f$; a proof of Theorem \rf{Te0.2} can be found in [Ar] p.~194. 

\sm The problem with Poincar\'e-Dulac normal forms is that, usually, they are not unique. In particular, one may wonder whether it could be possible to have such a normal form including finitely many resonant monomials only. This is indeed the case (see, e.g., Reich [Re]) when $\d f_O$ belongs to the so-called {\it Poincar\'e domain}, that is when $\d f_O$ is invertible and $O$ is either {\it attracting}, i.e., all the eigenvalues of~$\d f_O$ have modulus less than~$1$, or {\it repelling}, i.e., all the eigenvalues of~$\d f_O$ have modulus greater than~$1$ 
(when $\d f_O$ is still invertible but does not belong to the Poincar\'e domain, we shall say that it belongs to the {\it Siegel domain}). 

\sm Even without resonances, the holomorphic linearization is not guaranteed. The best positive result is due to Brjuno [Brj].
To describe Brjuno's result, let us introduce the
following definitions:

\sm\defin{De1.3.7} For $\lambda_1,\ldots,\lambda_n\in\C$ and $m\ge 2$ set
$$
\omega_{\lambda_1,\ldots,\lambda_n}(m)=\min_{2\le |Q|\le m\atop 1\le j\le n} |\lambda^Q-\lambda_j|.
\tag{eqcB}
$$
If $\lambda_1,\ldots,\lambda_n$ are the eigenvalues of~$\d f_O$, we shall
write~$\omega_f(m)$ for~$\omega_{\lambda_1,\ldots,\lambda_n}(m)$.

\sm It is clear that $\omega_f(m)\ne 0$ for all~$m\ge 2$ if and only if there are no
resonances. It is also not difficult to prove that if $f$ belongs to the
Siegel domain then
$$
\lim_{m\to+\infty}\omega_f(m)=0\;,
$$
which is the reason why, even without resonances, the formal linearization might
be diverging.

\sm\defin{De1.3.7bis} Let~$n\ge2$ and let~$\lambda_1, \dots, \lambda_n\in\C^*$ be not necessarily distinct. We say that~$\lambda$ {\it satisfies the Brjuno condition} if there exists a strictly increasing sequence of integers~$\{p_\nu\}_{\nu_\ge 0}$ with~$p_0=1$ such that
$$
\sum_{\nu\ge 0} {1\over p_\nu} \log{1\over\omega_{\lambda_1,\ldots,\lambda_n}(p_{\nu+1})}< \io.
\tag{eqcdue}
$$

Brjuno proved the following.

\sm\thm{Brjunohyp}{Theorem} (Brjuno, 1971 [Brj]) {\sl Let $f$ be a germ of biholomorphism of $\C^n$ fixing the origin, such that $\d f_O$ is diagonalizable. Assume moreover that the spectrum of $\d f_O$ has no resonances and satisfies the Brjuno condition.
%$$
%\sum_{k=0}^{+\infty}{1\over 2^k}\log{1\over\omega_f(2^{k+1})}<+\infty\;.
%\tag{eqcdue}
%$$
Then $f$ is holomorphically linearizable.}

%\sm\thm{Te0.3}{Theorem} (Poincar\'e, 1893 [Po]) {\sl Let~$f$ be a germ of holomorphic diffeomorphism of~$\C^n$ with an attracting or repelling fixed point. Then~$f$ is holomorphically linearizable if and only if it is formally linearizable. In particular, if there are no resonances then~$f$ is holomorphically linearizable.}

%\sm When~$O$ is not attracting or repelling, even without resonances, the formal linearization might diverge. In [Ra2] we found, under certain arithmetic conditions on the eigenvalues and some restrictions on the resonances, a necessary and sufficient condition for holomorphic linearization in presence of resonances, that in fact has as corollaries most of the known linearization results. In [Ra3] we found that, given~$m\ge 2$ germs~$f_1, \dots, f_m$ of biholomorphisms of~$\C^n$, fixing the origin, with~$(\d f_1)_O$ diagonalizable and such that~$f_1$ commutes with~$f_h$ for any~$h=2,\dots, m$, under certain arithmetic conditions on the eigenvalues of~$(\d f_1)_O$ and some restrictions on their resonances,~$f_1, \dots, f_m$ are simultaneously holomorphically linearizable if and only if there exists a particular complex manifold invariant under~$f_1, \dots, f_m$. 

\sm In the resonant case, one can still find formally linearizable germs, (see for example [R1] and [R2]), so two natural questions arise:

\sm
\item{(Q1)} {\it How many Poincar\'e-Dulac formal normal forms does a formally linearizable germ have?}

\sm
\item{(Q2)} {\it Is it possible to find arithmetic conditions on the eigenvalues of the spectrum of $\d f_O$ ensuring holomorphic linearizability of formally linearizable germs?}
\sm

R\"ussmann gave answers to both questions in [R\"u1], an I.H.E.S. preprint which is no longer available, and that was finally published in [R\"u2]. The answer to the first question is the following (the statement is slightly different from the original one presented in [R\"u2] but perfectly equivalent):

\sm\thm{TeFormalLin}{Theorem} (R\"ussmann, 2002 [R\"u2]) {\sl Let $f$ be a germ of biholomorphism of $\C^n$ fixing the origin. If $f$ is formally linearizable, then the linear form is its unique Poincar\'e-Dulac normal form.} 

\sm To answer to the second question, R\"ussmann introduced the following condition, that we shall call R\"ussmann condition.

\sm\defin{DeRussmann} Let~$n\ge2$ and let~$\lambda_1, \dots, \lambda_n\in\C^*$ be not necessarily distinct. We say that~$\lambda=(\lambda_1, \dots, \lambda_n)$ {\it satisfies the R\"ussmann condition} if there exists a function $\Omega\colon\N\to\R$ such that:{\parindent=35pt 
\item{(i)} $k\le \Omega(k)\le \Omega(k+1)$ for all $k\in\N$,
\sm\item{(ii)} $\sum\limits_{k\ge 1} {1\over k^2}\log\Omega(k)\le + \io$, and
\sm\item{(iii)} $|\lambda^Q -\lambda_j|\ge {1\over \Omega(|Q|)}$  for all $j=1,\dots n$ and for each multi-index $Q\in\N$ with $|Q|\ge 2$ not giving a resonance relative to $j$.\sm}

\sm R\"ussmann proved the following generalization of Brjuno's Theorem \rf{Brjunohyp}  (the statement is slightly different from the original one presented in [R\"u2] but perfectly equivalent).

\sm\thm{TeRussmann}{Theorem} (R\"ussmann, 2002 [R\"u2]) {\sl Let $f$ be a germ of biholomorphism of $\C^n$ fixing the origin and such that $\d f_O$ is diagonalizable. If $f$ is formally linearizable and the spectrum of $\d f_O$ satisfies the R\"ussmann condition, then $f$ is holomorphically linearizable.}

\sm We refer to [R\"u2] for the original proof and we limit ourselves to briefly recall here the main ideas. To prove these results, R\"ussmann first studies the process of Poincar\'e-Dulac formal normalization using a  functional iterative approach, without assuming anything on the diagonalizability of $\d f_O$. With this functional technique he proves Theorem \rf{TeFormalLin}; then he constructs a formal iteration process converging to a zero of the operator $\ca F(\phe) := f\circ \phe - \phe\circ \Lambda$ (where $\Lambda$ is the linear part of $f$), and, assuming $\Lambda$ diagonal, he gives estimates for each iteration step, proving that, under what we called the R\"ussmann condition, the process converges to a holomorphic linearization.

\me In this paper, we shall first present a direct proof of Theorem \rf{TeFormalLin} using power series expansions. Then we shall give a direct proof, using explicit computations with power series expansions and then proving convergence via majorant series, of an analogue of Theorem \rf{TeRussmann} under the following slightly different assumption, which is the natural generalization to the resonant case of the condition introduced by Brjuno. %Definitions \rf{De1.3.7} and \rf{De1.3.7bis} 

\sm\defin{De1.0} Let~$n\ge2$ and let~$\lambda_1, \dots, \lambda_n\in\C^*$ be not necessarily distinct. For $m\ge 2$ set
$$
\widetilde\omega_{\lambda_1,\ldots,\lambda_n}(m)=\min_{2\le |Q|\le m\atop Q\not\in {\res_j(\lambda)}} \min_{1\le j\le n} |\lambda^Q - \lambda_j|,
$$
where~$\res_j(\lambda)$ is the set of multi-indices~$Q\in\N^n$, with $|Q|\ge 2$, giving a resonance relation for~$\lambda =(\lambda_1, \dots, \lambda_n)$ relative to~$1\le j\le n$, i.e.,~$\lambda^Q-\lambda_j=0$.
If $\lambda_1,\ldots,\lambda_n$ are the eigenvalues of~$\d f_O$, we shall
write~$\widetilde\omega_f(m)$ for~$\widetilde\omega_{\lambda_1,\ldots,\lambda_n}(m)$.

%\sm{\bf inserire storia sulle risonanze...}

\sm\defin{De1.0bis} Let~$n\ge2$ and let~$\lambda=(\lambda_1, \dots, \lambda_n)\in(\C^*)^n$. We say that~$\lambda$ {\it satisfies the reduced Brjuno condition} if there exists a strictly increasing sequence of integers~$\{p_\nu\}_{\nu_\ge 0}$ with~$p_0=1$ such that
$$
\sum_{\nu\ge 0} {1\over p_\nu} \log{1\over \widetilde\omega_{\lambda_1,\ldots,\lambda_n}(p_{\nu+1})}< \io.
$$

We shall then prove:

\sm\thm{TeHolLin}{Theorem} {\sl Let $f$ be a germ of biholomorphism of $\C^n$ fixing the origin and such that $\d f_O$ is diagonalizable. If $f$ is formally linearizable and the spectrum of $\d f_O$ satisfies the reduced Brjuno condition, then $f$ is holomorphically linearizable.}

\sm We shall also show that R\"ussmann condition implies the reduced Brjuno condition and so our result implies Theorem \rf{TeRussmann}. The converse is known to be true in dimension $1$, as proved by R\"ussmann in [R\"u2], but is not known in higher dimension.

\me The structure of this paper is as follows. In the next section we shall discuss properties of formally linearizable germs, and we shall give our direct proof of Theorem \rf{TeFormalLin}. In section $3$ we shall prove Theorem \rf{TeHolLin} using majorant series. In the last section we shall discuss relations between R\"ussmann condition and the reduced Brjuno condition.

\me\no{\bf Acknowledgments.} I would like to thank Marco Abate for helpful comments on a draft of this work.

%%%%%%%%%%%%%%%%%%%%%%%%%%%%%%%%%%%%%%%%%%%%%%%%%%%%%%%%%%%%%%%%%%%%%%%%%%%%%%%%
\sect Formally linearizable germs
%%%%%%%%%%%%%%%%%%%%%%%%%%%%%%%%%%%%%%%%%%%%%%%%%%%%%%%%%%%%%%%%%%%%%%%%%%%%%%%%

In general, a germ $f$ can have several Poincar\'e-Dulac formal normal forms; however, we can say something on the shape of the formal conjugations between them. We have in fact the following result. 

\sm\thm{PrFormeConiugate}{Proposition} {\sl Let $f$ and $g$ be two germs of biholomorphism of $\C^n$ fixing the origin, with the same linear part $\Lambda$ and in Poincar\'e-Dulac normal form. If there exists a formal transformation $\phe$ of $\C^n$, with no constant term and tangent to the identity, conjugating $f$ and $g$, then $\phe$ contains only monomials that are resonant with respect to the eigenvalues of $\Lambda$.} 

\sm\proof Since $f$ and $g$ are in Poincar\'e-Dulac normal form, $\Lambda$ is in Jordan normal form. Let $\lambda_1,\dots,\lambda_n$ be the eigenvalues of $\Lambda$. We shall prove that a formal solution $\phe= I + \widehat\phe$ of 
$$
f\circ\phe= \phe\circ g
\tag{eqconiugio}
$$ 
contains only monomials that are resonant with respect to $\lambda_1,\dots, \lambda_n$. Using the standard multi-index notation, for each $j\in\{1,\dots, n\}$ we can write
$$
f_j(z)= \lambda_j z_j + \eps_j z_{j-1}+  z_j f_j^{\rm res}(z) = \lambda_j z_j + \eps_j z_{j-1}+  z_j \sum_{Q\in N_j \atop \lambda^Q=1} f_{Q,j} z^Q,
$$
$$
g_j(z)= \lambda_j z_j + \eps_j z_{j-1}+  z_j g_j^{\rm res}(z) =\lambda_j z_j + \eps_j z_{j-1} + z_j \sum_{Q\in N_j \atop \lambda^Q=1} g_{Q,j} z^Q,
$$
and
$$
\phe_j(z)= z_j\left( 1 + \phe_j^{\rm res}(z) + \phe_j^{\ne{\rm res}}(z)\right)= z_j + z_j \sum_{Q\in N_j \atop \lambda^Q=1} \phe_{Q,j} z^Q + z_j \sum_{Q\in N_j \atop \lambda^Q\ne1} \phe_{Q,j} z^Q,
$$
where 
$$
N_j:=\{Q\in \Z^n \mid |Q|\ge 1, q_j\ge -1, q_h\ge 0~\hbox{for all}~h\ne j\},
$$
and $\eps_j\in\{0,1\}$ can be non-zero only if $\lambda_j = \lambda_{j-1}$. With these notations, the left-hand side of the $j$-th coordinate of \rf{eqconiugio} becomes
$$
\eqalign{
(f\circ \phe)_j(z) 
&= \lambda_j\phe_j(z) + \eps_j\phe_{j-1}(z) + \phe_j(z) \sum_{Q\in N_j\atop \lambda^Q=1} f_{Q,j} \prod_{k=1}^n \phe_k(z)^{q_k}\cr
&= \lambda_j z_j\left( 1 + \phe_j^{\rm res}(z) + \phe_j^{\ne{\rm res}}(z)\right)\cr
&\quad + \eps_j z_{j-1}\left( 1 + \phe_{j-1}^{\rm res}(z) + \phe_{j-1}^{\ne{\rm res}}(z)\right)\cr
&\quad + z_j\left( 1 + \phe_j^{\rm res}(z) + \phe_j^{\ne{\rm res}}(z)\right)\!\!\!\sum_{Q\in N_j\atop \lambda^Q=\lambda_j} \!\!f_{Q,j} z^Q \prod_{k=1}^n \left( 1 + \phe_k^{\rm res}(z) + \phe_k^{\ne{\rm res}}(z)\right)^{q_k}\!\!,
}\tag{eqconiugiosinistra}
$$
while the $j$-th coordinate of the right-hand side of \rf{eqconiugio} becomes 
$$
\eqalign{
(\phe\circ g)_j(z) 
&= g_j(z) + g_j(z) \sum_{Q\in N_j\atop \lambda^Q=1} \phe_{Q,j} \prod_{k=1}^n g_k(z)^{q_k} + g_j(z)\sum_{Q\in N_j\atop \lambda^Q\ne1} \phe_{Q,j} \prod_{k=1}^n g_k(z)^{q_k}\cr
&= \lambda_j z_j + \eps_j z_{j-1}+  z_j g_j^{\rm res}(z) \cr
&\quad+ \left(\lambda_j z_j + \eps_j z_{j-1}+  z_j g_j^{\rm res}(z)\right) \!\!\!\sum_{Q\in N_j\atop \lambda^Q=1} \!\!\phe_{Q,j}z^Q \prod_{k=1}^n \left(\lambda_k + \eps_k {z_{k-1}\over z_k} +  g_k^{\rm res}(z)\right) ^{q_k}\cr
&\quad+ \left(\lambda_j z_j + \eps_j z_{j-1}+  z_j g_j^{\rm res}(z)\right) \!\!\!\sum_{Q\in N_j\atop \lambda^Q \ne 1}\!\!\phe_{Q,j}z^Q  \prod_{k=1}^n \left(\lambda_k + \eps_k {z_{k-1}\over z_k} +  g_k^{\rm res}(z)\right) ^{q_k}\!\!\!\!.
}\tag{eqconiugiodestra}
$$
%$$
%\eqalign{
%(f\circ \phe)_j(z) 
%&= z_j\!\!\left(\!\!1 \! + \!\!\! \sum_{P\in N_j \atop \lambda^P=1} \!\phe_{P,j} z^P\!\! + \!\!\!\sum_{P\in N_j \atop \lambda^P\ne1} \!\phe_{P,j} z^P\!\!\right) \!\!\!\left(\!\lambda_j \! + \!\!\sum_{Q\in N_j \atop \lambda^Q=1} f_{Q,j} z^Q \!\prod_{k=1}^n\!\!\left(\!\!1 \!+ \!\!\!\sum_{P\in N_k \atop \lambda^P=1} \!\phe_{P,k} z^P \!\!+\!\!\! \sum_{P\in N_k \atop \lambda^P\ne1} \!\phe_{P,k} z^P\!\!\right)^{\!\!\!\!q_k}\!\right)\cr
%&= z_j\left(\lambda_j + \sum_{P\in N_j \atop \lambda^P=1} g_{P,j} z^P\right) \left(1 + \sum_{Q\in N_j \atop \lambda^Q=1} \phe_{Q,j} z^Q\prod_{k=1}^n\left(\lambda_j + \sum_{P\in N_k \atop \lambda^P=1} g_{P,k} z^P\right)^{q_k}\right .\cr
%&\qquad\qquad\qquad\qquad\qquad\qquad\qquad\qquad+\left . \sum_{Q\in N_j \atop \lambda^Q\ne1} \phe_{Q,j} z^Q\prod_{k=1}^n\left(\lambda_j + \sum_{P\in N_k \atop \lambda^P=1} g_{P,k} z^P\right)^{q_k}\right)\cr
%&=(\phe\circ g )_j(z). 
%}
%$$
Furthermore, notice that if $P$ and $Q$ are two multi-indices such that $\lambda^P=\lambda^Q=1$, then we have $\lambda^{\alpha P + \beta Q} = 1$ for every $\alpha, \beta\in \Z$. 

We want to prove that $\phe_{Q,j}=0$ for each multi-index $Q\in N_j$ such that $\lambda^Q\ne 1$.
Let us assume by contradiction that this is not true, and let $\widetilde Q$ be the first (with respect to the lexicographic order) multi-index in $N:=\bigcup_{j=1}^n N_j$ so that $\lambda^{\widetilde Q}\ne 1$ and $\phe_{\widetilde Q,j}\ne 0$. Let $j$ be the minimal in $\{1,\dots, n\}$ such that $\widetilde Q\in N_j$, and let us compute the coefficient of the monomial $z^{\widetilde Q+e_j}$ in \rf{eqconiugiosinistra} and \rf{eqconiugiodestra}. In \rf{eqconiugiosinistra} we only have $\lambda_j\phe_{\widetilde Q, j}$ because, since~$f-\Lambda$ is of second order and resonant,
other contributions could come only from coefficients~$\psi_{P,k}$ with~$|P| < |\widetilde Q|$ and $\lambda^P\ne1$, but there are no such coefficients thanks to the minimality of $\widetilde Q$ and $j$. In \rf{eqconiugiodestra} we can argue analogously, but we have also to take care of the monomials divisible by $\eps_k^h(z_{k-1}/z_k)^h z^P$, with $\lambda^P=1$; in this last case, if $\eps_k\ne0$, we obtain a multi-index $P-he_k+ h e_{k-1}$, and again $\lambda^{P-he_k+ h e_{k-1}}=1$ because $\lambda_k=\lambda_{k+1}$. Then in \rf{eqconiugiodestra} we only have $\lambda^{\widetilde Q+e_j}\phe_{\widetilde Q, j}$. Hence, we have
$$
(\lambda^{\widetilde Q + e_j} -\lambda_j)\phe_{\widetilde Q, j} = 0,
$$
yielding 
$$
\phe_{\widetilde Q, j} = 0,
$$
because $\lambda^{\widetilde Q } \ne 1$ and $\lambda_j\ne 0$, and contradicting the hypothesis. \qed

\sm\thm{Re1.3.25}{Remark} It is clear from the proof that Proposition \rf{PrFormeConiugate} holds also in the formal category, i.e., for $f, g\in\C_O[\![z_1,\dots,z_n]\!]$ formal power series without constant terms in Poincar\'e-Dulac normal form.

\sm We can now give a direct proof of Theorem \rf{TeFormalLin}, i.e., that when a germ is formally linearizable, then the linear form is its unique Poincar\'e-Dulac normal form.

\sm\thm{TeFormalLin}{Theorem} {\sl Let $f$ be a germ of biholomorphism of $\C^n$ fixing the origin. If $f$ is formally linearizable, then the linear form is its unique Poincar\'e-Dulac normal form.} 
%{\sl Let $f$ be a germ of biholomorphism of $\C^n$ fixing the origin. If $f$ is formally linearizable, and $f$ is formally conjugated to a formal Poincar\'e-Dulac normal form $g$, then $g$ is linear.} 

\sm\proof Let $\Lambda$ be the linear part of $f$. Up to linear conjugacy, we may assume that $\Lambda$ is in Jordan normal form. If the eigenvalues $\lambda_1,\dots,\lambda_n$ of $\Lambda$ have no resonances, then there is nothing to prove. Let us then assume that we have resonances, and let us assume by contradiction that there is another Poincar\'e-Dulac formal normal form $g\not\equiv \Lambda$ associated to $f$. Since $f$ is formally linearizable and it is formally conjugated to $g$, also $g$ is formally linearizable. Thanks to Proposition \rf{PrFormeConiugate}, any formal linearization $\psi$ of $g$ tangent to the identity contains only monomials resonant with respect to $\lambda_1,\dots,\lambda_n$; hence, writing $g = \Lambda + g^{\rm res}$ and $\psi = I + \psi^{\rm res}$, the conjugacy equation $g\circ \psi = \psi\circ \Lambda$ becomes
$$
\eqalign{\Lambda + \Lambda\psi^{\rm res} + g^{\rm res}\circ (I + \psi^{\rm res}) 
&=(\Lambda + g^{\rm res})\circ (I + \psi^{\rm res})\cr 
&=(I +\psi^{\rm res})\circ\Lambda\cr
&=\Lambda + \psi^{\rm res}\circ\Lambda\cr
&=\Lambda + \Lambda \psi^{\rm res},}
$$ 
because $\psi^{\rm res}\circ\Lambda=\Lambda\psi^{\rm res}$.
Hence there must be
$$
g^{\rm res}\circ \psi\equiv 0,
%\tag{eqFormalLin}
$$
and composing on the right with $\psi^{-1}$ we get $g^{\rm res}\equiv 0$. \qed

\sm\thm{Re1.4.2.1}{Remark} As a consequence of the previous result, we get that {\sl any formal normalization given by the Poincar\'e-Dulac procedure applied to a formally linerizable germ $f$ is indeed a formal linearization of the germ.} In particular, we have {\sl uniqueness of the Poincar\'e-Dulac normal form (which is linear and hence holomorphic), but not of the formal linearizations}. Hence {\sl a formally linearizable germ $f$ is formally linearizable via a formal transformation $\phe= \id + \widehat\phe$ containing only non-resonant monomials.} In fact, thanks to the standard proof of Poincar\'e-Dulac Theorem (see [R3] Theorem $1.3.25$), we can consider the formal normalization obtained with the Poincar\'e-Dulac procedure and imposing $\phe_{Q,j} = 0$ for all $Q$ and $j$ such that $\lambda^Q=\lambda_j$; and this formal transformation $\phe$, by Theorem~\rf{TeFormalLin}, conjugates $f$ to its linear part.

%%%%%%%%%%%%%%%%%%%%%%%%%%%%%%%%%%%%%%%%%%%%%%%%%%%%%%%%%%%%%%%%%%%%%%%%%%%%%%%%
\sect Convergence under the reduced Brjuno condition
%%%%%%%%%%%%%%%%%%%%%%%%%%%%%%%%%%%%%%%%%%%%%%%%%%%%%%%%%%%%%%%%%%%%%%%%%%%%%%%%

Now we have all the ingredients needed to prove Theorem \rf{TeHolLin}.

\sm\thm{TeHolLin}{Theorem} {\sl Let $f$ be a germ of biholomorphism of $\C^n$ fixing the origin and such that $\d f_O$ is diagonalizable. If $f$ is formally linearizable and the spectrum of $\d f_O$ satisfies the reduced Brjuno condition, then $f$ is holomorphically linearizable.}

\sm\proof Up to linear changes of the coordinates, we may assume that the linear part $\Lambda$ of $f$ is diagonal, i.e., $\Lambda = \diag(\lambda_1,\dots, \lambda_n)$. From the conjugacy equation
$$
f\circ \phe=\phe\circ \Lambda,
\tag{eq1.1}
$$
writing $f(z)=\Lambda z + \sum_{|L|\ge 2} f_L z^L$, and $\phe(w) = w + \sum_{|Q|\ge 2} \phe_Q w^Q$, where $f_L$ and $\phe_Q$ belong to $\C^n$, we have that coefficients of $\phe$ have to verify   
$$
\sum_{|Q|\ge 2} A_Q \phe_Q w^Q = \sum_{|L|\ge 2}f_L \left(\sum_{|M|\ge1}\phe_M w^M \right)^L, 
\tag{eq1.3}
$$
where
$$
A_Q=\lambda^Q I_n-\Lambda.
$$
The matrices~$A_Q$ are not invertible only when $Q\in\bigcup_{j=1}^n \res_j(\lambda)$, but, thanks Remark \rf{Re1.4.2.1}, we can set $\phe_{Q,j}=0$ for all $Q\in\res_j(\lambda)$; hence we just have to consider $Q\not\in\bigcap_{j=1}^n \res_j(\lambda)$, and, to prove the convergence of the formal conjugation~$\phe$ in a neighbourhood of the origin, it suffices to show that
$$
\sup_Q{1\over |Q|} \log\|\phe_Q\|<\io,
\tag{eq6}
$$
for $|Q|\ge 2$ and $Q\not\in\cap_{j=1}^n \res_j(\lambda)$.

Since~$f$ is holomorphic in a neighbourhood of the origin, there exists a positive number~$\rho$ such that~$\|f_L\|\le \rho^{|L|}$ for~$|L|\ge 2$. The functional equation \rf{eq1.1} remains valid under the linear change of coordinates~$f(z)\mapsto \sigma f(z/\sigma)$,~$\phe(w)\mapsto \sigma\phe(w/\sigma)$ with~$\sigma=\max\{1, \rho^2\}$. Therefore we may assume that 
$$
\forevery{|L|\ge 2}{\|f_L\|\le 1.}
$$
It follows from \rf{eq1.3} that for any multi-index $Q\in\N^n\setminus\bigcap_{j=1}^n \res_j(\lambda)$ with $|Q|\ge 2$ we have 
$$
\|\phe_Q\|\le \eps_Q^{-1} \sum_{Q_1+\cdots +Q_\nu = Q \atop \nu\ge 2} \|\phe_{Q_1}\|\cdots \|\phe_{Q_\nu}\|,
\tag{eq7}
$$
where
$$
\eps_Q = \min_{1\le j \le n\atop Q \not\in\res_j(\lambda)} |\lambda^Q - \lambda_j|.
$$
We can define, inductively, for~$m \ge 2$ 
$$
\alpha_m= \sum_{m_1+\cdots + m_\nu =j \atop \nu \ge 2} \alpha_{m_1} \cdots \alpha_{m_\nu},
$$
and
$$
\delta_Q = \eps_Q^{-1}\max_{Q_1+\cdots + Q_\nu =Q\atop \nu\ge 2} \delta_{Q_1}\cdots\delta_{Q_\nu},
$$
for $Q\in\N^n\setminus \bigcap_{j=1}^n \res_j(\lambda)$ with~$|Q|\ge 2$,
with~$\alpha_1 =1$ and~$\delta_E= 1$, where~$E$ is any integer vector with~$|E|=1$. Then, by induction, we have that
$$
{\|\phe_Q\|\le \alpha_{|Q|}\delta_Q,}
$$
for every $Q\in\N^n\setminus \bigcap_{j=1}^n \res_j(\lambda)$ with~$|Q|\ge 2$. Therefore, to establish \rf{eq6} it suffices to prove analogous estimates for~$\alpha_m$ and~$\delta_Q$.

\sm It is easy to estimate~$\alpha_m$. Let~$\alpha= \sum_{m \ge 1}\alpha_m t^m$. We have
$$
\eqalign{\alpha - t &= \sum_{m\ge 2} \alpha_m t^m \cr
					&= \sum_{m\ge 2}\left(\sum_{h\ge 1}\alpha_h t^h\right)^m\cr
					&= {\alpha^2 \over 1- \alpha}.}
$$
This equation has a unique holomorphic solution vanishing at zero
$$
\alpha= {t+1 \over 4} \left(1 - \sqrt{1-{8t\over(1+t)^2}}\right),
$$
defined for~$|t|$ small enough. Hence,
$$
\sup_m {1\over m}\log \alpha_m < \io,
$$
as we want.

\sm To estimate~$\delta_Q$ we have to take care of small divisors.
First of all, for each multi-index~$Q\not\in \bigcap_{j=1}^n \res_j(\lambda)$ with~$|Q|\ge 2$ we can associate to~$\delta_Q$ a decomposition of the form
$$
\delta_Q= \eps_{L_0}^{-1}\eps_{L_1}^{-1}\cdots\eps_{L_p}^{-1},\tag{eqdelta}
$$
where~$L_0=Q$,~$|Q|>|L_1|\ge\cdots\ge|L_p|\ge2$ and~$L_j\not\in \bigcap_{j=1}^n \res_j(\lambda)$ for all~$j=1, \dots, p$ and~$p\ge 1$. 
In fact, we choose a decomposition~$Q=Q_1+\cdots+Q_\nu$ such that the maximum in the expression of~$\delta_Q$ is achieved; obviously,~$Q_j$ does not belong to~$\bigcap_{j=1}^n \res_j(\lambda)$ for all~$j=1,\dots, \nu$. We can then express~$\delta_Q$ in terms of~$\eps_{Q_j}^{-1}$ and~$\delta_{Q'_j}$  with~$|Q'_j|<|Q_j|$. Carrying on this process, we eventually arrive at a decomposition of the form \rf{eqdelta}. Furthermore, for each multi-index $ Q\not\in \bigcap_{j=1}^n \res_j(\lambda)$ with $|Q|\ge 2$, we can choose an index $i_Q$ so that
$$
\eps_Q = |\lambda^Q - \lambda_{i_Q}|.
$$ 

\sm The rest of the proof follows closely the proof of Theorem $5.1$ in [R1]. For the benefit of the reader, we report it here.

For $m\ge 2$ and $1\le j\le n$, we can define
$$
N^j_m(Q)
$$
to be the number of factors~$\eps_{L}^{-1}$ in the expression \rf{eqdelta} of~$\delta_Q$, satisfying
$$
\eps_{L}<\theta\,\widetilde\omega_f(m),~~\hbox{and}~~i_L = j,
$$ 
where~$\widetilde\omega_f(m)$ is defined in Definition \rf{De1.0}, and in this notation can be expressed as
$$
\widetilde\omega_f(m)= \min_{2\le|Q|\le m \atop Q\not\in {\cap}_{j=1}^n \res_j(\lambda)}\eps_Q, 
$$
and~$\theta$ is the positive real number satisfying
$$
4\theta=\min_{1\le h\le n}|\lambda_h| \le 1.
$$
The last inequality can always be satisfied by replacing~$f$ by~$f^{-1}$ if necessary. Moreover we also have~$\widetilde\omega_f(m)\le 2$.

Notice that~$\widetilde\omega_f(m)$ is non-increasing with respect to~$m$ and under our assumptions~$\widetilde\omega_f(m)$ tends to zero as~$m$ goes to infinity. The following is the key estimate. 

\sm\thm{Le1.1}{Lemma} {\sl For~$m\ge2$,~$1\le j\le n$ and~$Q\not\in\bigcap_{j=1}^n \res_j(\lambda)$, we have
$$
N^j_m(Q)\le \cases{0, &\hbox{if}~~$|Q|\le m$,\cr\noalign{\sm}
					\displaystyle {2|Q|\over m}-1,&\hbox{if}~~$|Q|> m$.}
$$}
\sm\proof The proof is done by induction on $|Q|$. Since we fix~$m$ and~$j$ throughout the proof, we write~$N$ instead of~$N^j_m$.

For~$|Q|\le m$, 
$$
\eps_Q\ge\widetilde\omega_f(|Q|)\ge \widetilde\omega_f(m) > \theta\, \widetilde\omega_f(m),
$$
hence~$N(Q)=0$.

Assume now that~$|Q|>m$. Then~$2|Q|/m -1 \ge 1$. Write
$$
\delta_Q= \eps_Q^{-1}
\delta_{Q_1}\cdots \delta_{Q_\nu}, \quad Q=Q_1 + \cdots + Q_\nu, \quad \nu\ge2,
$$
with~$|Q|>|Q_1|\ge \cdots\ge|Q_\nu|$; note that~$Q - Q_1$ does not belong to~$\bigcap_{j=1}^n \res_j(\lambda)$, otherwise the other~$Q_h$'s would be in $\bigcap_{j=1}^n \res_j(\lambda)$. We have to consider the following different cases.

\sm{\it Case 1:}~$\eps_Q \ge \theta\,\widetilde\omega_f(m)$ and~$i_Q$ arbitrary, or~$\eps_Q < \theta\,\widetilde\omega_f(m)$ and~$i_Q\ne j$. Then
$$
N(Q) = N(Q_1) + \cdots + N(Q_\nu),
$$
and applying the induction hypotheses to each term we get~$N(Q)\le (2|Q|/m) - 1$.

\sm{\it Case 2:}~$\eps_Q<\theta\,\widetilde\omega_f(m)$ and~$i_Q=j$. Then
$$
N(Q) = 1 + N(Q_1) + \cdots + N(Q_\nu),
$$
and there are three different subcases.

{\it Case 2.1:}~$|Q_1|\le m$. Then
$$
N(Q) = 1 < {2|Q|\over m} -1,
$$
as we want.

{\it Case 2.2:}~$|Q_1|\ge|Q_2|>m$. Then there is~$\nu'$ such that~$2\le\nu'\le\nu$ and~$|Q_{\nu'}|> m\ge |Q_{\nu'+1}|$, and we have
$$
N(Q) = 1 + N(Q_1) + \cdots + N(Q_{\nu'})\le 1 + {2|Q|\over m} - \nu' \le  {2|Q|\over m} -1.
$$

{\it Case 2.3:}~$|Q_1|>m\ge |Q_2|$. Then
$$
N(Q)= 1 + N(Q_1),
$$
and there are again three different subcases.

{\it Case 2.3.1:}~$i_{Q_1}\ne j$. Then~$N(Q_1) = 0$ and we are done. 

{\it Case 2.3.2:}~$|Q_1|\le|Q|-m$ and~$i_{Q_1} = j$. Then
$$
N(Q) \le 1 + 2\,{|Q|-m\over m} -1 < {2|Q|\over m} -1.
$$

{\it Case 2.3.3:}~$|Q_1|>|Q|-m$ and~$i_{Q_1} = j$. The crucial remark is that~$\eps_{Q_1}^{-1}$ gives no contribute to~$N(Q_1)$, as shown in the next lemma.

\sm\thm{Le1.2}{Lemma} {\sl If~$Q>Q_1$ with respect to the lexicographic order,~$Q$,~$Q_1$ and~$Q-Q_1$ are not in~$\bigcap_{j=1}^n \res_j(\lambda)$,~$i_Q= i_{Q_1} =j$ and
$$
\eps_Q < \theta\, \widetilde\omega_f(m) \quad\hbox{and}\quad \eps_{Q_1}<\theta\, \widetilde\omega_f(m),
$$
then~$|Q-Q_1| = |Q| - |Q_1| \ge m$.}

\sm\proof Before we proceed with the proof, notice that the equality~$|Q-Q_1| = |Q| - |Q_1|$ is obvious since~$Q>Q_1$.

Since we are supposing~$\eps_{Q_1} = |\lambda^{Q_1}-\lambda_j| <\theta\,\widetilde \omega_f(m)$, we have
$$
\eqalign{|\lambda^{Q_1}| &>|\lambda_j|-\theta\,\widetilde\omega_f(m)\cr
								&\ge 4\theta - 2\theta = 2\theta.
}$$
Let us suppose by contradiction~$|Q-Q_1| = |Q| - |Q_1| < m$. By assumption, it follows that
$$
\eqalign{2\theta\,\widetilde\omega_f(m) &> \eps_Q+ \eps_{Q_1}\cr
							&= |\lambda^Q -\lambda_{j}| + |\lambda^{Q_1} -\lambda_{j}|\cr
							&\ge |\lambda^Q - \lambda^{Q_1}|\cr
							&\ge |\lambda^{Q_1}| \,|\lambda^{Q-Q_1} - 1|\cr
							&\ge 2\theta\,\widetilde \omega_f(|Q-Q_1|+1)\cr
							&\ge 2\theta\, \widetilde\omega_f(m),
}$$
which is impossible. \qed

\sm Using Lemma~\rf{Le1.2}, case~$1$ applies to~$\delta_{Q_1}$ and we have
$$
N(Q) = 1 + N(Q_{1_1}) + \cdots + N(Q_{1_{\nu_{1}}}),
$$
where~$|Q|> |Q_1|> |Q_{1_1}| \ge \cdots \ge |Q_{1_{\nu_{1}}}|$ and~$Q_1 =Q_{1_1} + \cdots +Q_{1_{\nu_{1}}}$. We can do the analysis of case~$2$ again for this decomposition, and we finish unless we run into case~$2.3.2$ again. However, this loop cannot happen more than~$m+1$ times and we have to finally run into a different case. This completes the induction and the proof of Lemma~\rf{Le1.1}. \qed

\sm Since the spectrum of~$\d f_O$ satisfies the reduced Brjuno condition, there exists a strictly increasing sequence~$\{p_\nu\}_{\nu\ge 0}$ of integers with~$p_0=1$ and such that
$$
\sum_{\nu\ge 0} {1\over p_\nu}\log{1\over\widetilde\omega_f(p_{\nu +1})} <\io.\tag{eq8}
$$
We have to estimate
$$
{1\over |Q|}\log\delta_Q = \sum_{j=0}^p {1\over |Q|} \log \eps_{L_j}^{-1}, \quad Q\not\in\bigcap_{j=1}^n \res_j(\lambda).
$$
By Lemma~\rf{Le1.1}, 
$$
\eqalign{{\rm card}\left\{0\le j\le p : \theta\, \widetilde\omega_f(p_{\nu +1}) \le \eps_{L_j} <\theta\, \widetilde\omega_f(p_\nu)\right\} &\le N_{p_\nu}^1(Q) + \cdots N_{p_\nu}^n(Q) \cr
							&\le {2n|Q|\over p_\nu}} 
$$
for~$\nu\ge 1$. It is also easy to see from the definition of~$\delta_Q$ that the number of factors~$\eps_{L_j}^{-1}$ is bounded by~$2|Q| - 1$. In particular,
$$
{\rm card}\left\{0\le j\le p : \theta\, \widetilde\omega_f(p_{1}) \le \eps_{L_j} \right\} \le 2n|Q| = {2n|Q|\over p_0}. 
$$
Then,
$$
\eqalign{{1\over |Q|} \log \delta_Q &\le 2n \sum_{\nu\ge 0} {1\over p_\nu} \log{1\over\theta\,\widetilde\omega_f(p_{\nu +1})} \cr
	&= 2n\left( \sum_{\nu \ge 0} {1\over p_\nu}\log{1\over\widetilde\omega_f(p_{\nu +1})} + \log{1\over\theta} \sum_{\nu \ge 0} {1\over p_\nu}\right).}\tag{eq9}
$$
Since~$\widetilde\omega_f(m)$ tends to zero monotonically as~$m$ goes to infinity, we can choose some~$\overline{m}$ such that~$1>\widetilde\omega_f(m)$ for all~$m>\overline{m}$, and we get
$$
\sum_{\nu\ge\nu_0} {1\over p_\nu} \le {1\over \log (1/\widetilde\omega_f(\overline{m}))} \sum_{\nu\ge \nu_0} {1\over p_\nu} \log{1\over\widetilde\omega_f(p_{\nu +1})},
$$
where~$\nu_0$ verifies the inequalities~$p_{\nu_0 -1}\le \overline{m} < p_{\nu_0}$. Thus both series in parentheses in \rf{eq9} converge thanks to \rf{eq8}. Therefore
$$
\sup_Q {1\over |Q|}\log \delta_Q <\io
$$ 
and this concludes the proof. \qed

\sm When there are no resonances, we obtain Brjuno's Theorem \rf{Brjunohyp}. 

\sm\thm{ReNotLin}{Remark} If the reduced Brjuno condition is not satisfied, then there are formally linearizable germs that are not holomorphically linearizable. A first example is the following: let us consider the following germ of biholomorphism $f$ of $(\C^2,O)$:
$$
\eqalign{&f_1(z,w)=\lambda z + z^2,\cr
  		 &f_2(z,w)= w,}\tag{ex1}
$$
with $\lambda = e^{2\pi i \theta}$, $\theta\in \R\setminus\Q$, not a Brjuno number. We are in presence of resonances because~$\res_1(\lambda, 1) = \{P\in\N^2\mid P = (1, p), p\ge 1\}$ and $\res_2(\lambda, 1) = \{P\in\N^2\mid P = (0, p), p\ge 2\}$. It is easy to prove that $f$ is formally linearizable, but not holomorphically linearizable, because otherwise the holomorphic function~$\lambda z + z^2$ would be holomorphically linearizable contradicting Yoccoz's result [Y].

\sm A more general example is the following:

\sm\thm{ExNotLin}{Example} Let $n\ge 2$, and let $\lambda_1,\dots,\lambda_s\in\C^*$, be $1\le s<n$ complex non-resonant numbers such that
$$
\limsup_{m\to+\infty}{1\over
m}\log{1\over\omega_{\lambda_1,\ldots,\lambda_s}(m)}=+\infty\;.
\tag{eqcremer}
$$
Then it is possible to find (see e.g. [R3] Theorem $1.5.1$) a germ $f$ of biholomorphism of $\C^s$ fixing the origin, with $\d f_O=\diag(\lambda_1,\ldots,\lambda_s)$, formally linearizable (since there are no resonances) but not holomorphically linearizable. It is also possible to find $\mu_1,\dots, \mu_r\in\C^*$, with $r=n-s$, such that the $n$-tuple $\lambda=(\lambda_1, \dots, \lambda_s,\mu_1,\dots, \mu_r)\in(\C^*)^n$ {\it has only level $s$ resonances} (see [R1], where this definition was first introduced, for details), i.e., for $1\le j\le s$ we have 
$$
\res_j(\lambda) = \{P\in \N^n \mid |P|\ge 2,\, p_l =\delta_{jl}~\hbox{for}~l=1,\dots, s,~~\hbox{and}~~\mu_1^{p_{s+1}}\cdots \mu_r^{p_n}=1\}, 
$$
where $\delta_{jl}$ is the Kroenecker's delta, and for $s+1\le h\le n$ we have 
$$
\res_h(\lambda)= \{P\in \N^n \mid |P|\ge2,~p_1=\cdots=p_s=0,~\mu_1^{p_{s+1}}\cdots \mu_r^{p_n}=\mu_{h-s}\}.
$$
Then any germ of biholomorphism $F$ of $\C^n$ fixing the origin of the form %f_j(z) = \lambda_j z_j + \widetilde f_j(z) 
$$
\casi{
F_j(z,w) = f_j(z) &\hbox{for}~$j=1,\dots, s$,\cr\noalign{\sm}
F_h(z,w) = \mu_{h-s} w_{h-s} + \widetilde F_h(z,w) &\hbox{for}~$h=s+1,\dots, n$, 
}
$$
with 
$$
{\rm ord}_z(\widetilde F_h)\ge 1,
$$
for $h=s+1, \dots, n$, where $(z,w) = (z_1,\dots, z_s, w_1, \dots w_r)$ are local coordinates of $\C^n$ at the origin, is formally linearizable (see Theorem $4.1$ of [R1]), but $\lambda= (\lambda_1, \dots, \lambda_s,\mu_1,\dots, \mu_r)$ does not satisfy the reduced Brjuno condition (because of \rf{eqcremer}) and $F$ is not holomorphically linearizable. In fact, if $F$ were holomorphically linearizable via a linearization $\Phi$, tangent to the identity, then $F\circ \Phi = \Phi \circ \diag(\lambda_1, \dots, \lambda_s,\mu_1,\dots, \mu_r)$. Hence, for each $1\le j\le s$, we would have
$$
\eqalign{
(F\circ \Phi)_j(z,w) 
&= \lambda_j \Phi_j(z,w) + \widetilde f_j(\Phi_1(z,w), \dots, \Phi_s(z,w))
\cr
&=(\Phi\circ \diag(\lambda_1, \dots, \lambda_s,\mu_1,\dots, \mu_r))_j(z,w)\cr
&=\Phi_j(\lambda_1 z_1, \dots, \lambda_s z_s,\mu_1 w_1,\dots, \mu_r w_r),
}
$$
yielding
$$
\eqalign{
(F\circ \Phi)_j(z,0) 
%&= \lambda_j \Phi_j(z,0) + \widetilde f_j(\Phi_1(z,0), \dots, \Phi_s(z,0))
%\cr
%&=(\Phi\circ \diag(\lambda_1, \dots, \lambda_s,\mu_1,\dots, \mu_r))_j(z,0)\cr
&=\Phi_j(\lambda_1 z_1, \dots, \lambda_s z_s,0,\dots, 0),
}
$$ 
and thus the holomorphic germ $\phe$ of $\C^s$ fixing the origin defined by $\phe_j(z) = \Phi_j(z,0)$ for $j=1,\dots, s$, would coincide with the unique formal linearization of $f$, that would then be convergent contradicting the hypotheses.

%%%%%%%%%%%%%%%%%%%%%%%%%%%%%%%%%%%%%%%%%%%%%%%%%%%%%%%%%%%%%%%%%%%%%%%%%%%%%%%%
\sect R\"ussmann condition vs. reduced Brjuno condition
%%%%%%%%%%%%%%%%%%%%%%%%%%%%%%%%%%%%%%%%%%%%%%%%%%%%%%%%%%%%%%%%%%%%%%%%%%%%%%%%

R\"ussmann proves that, in dimension $1$, his condition is equivalent to Brjuno condition (see Lemma~8.2 of [R\"u2]), and he also proves the following result.

%\bi{\bf and it also ``proves'' that its condition ``implies'' the Brjuno reduced condition? da controllare se \`e vero ;-), da rielaborare per ottenere la condizione di Brjuno visto che R\"ussmann lo fa usando la sua $\Omega$ e poi introdurre le varie condizioni di Brjuno}

\sm\thm{LeRussmann}{Lemma} (R\"ussmann, 2002 [R\"u2]) {\sl Let $\Omega\colon\N\to(0,+\io)$ be a monotone non decreasing function, and 
let $\{s_\nu\}$ be defined by $s_\nu := 2^{q+\nu}$, with $q\in\N$. Then
$$
\sum_{\nu\ge 0} {1\over s_\nu} \log \Omega(s_{\nu+1}) \le \sum_{k\ge 2^{q+1}} {1\over k^2} \log \Omega(k).
$$ 
}

%\sm\proof For each $a$, $b$ integers with $0< a< b$ we have
%$$
%{1\over a} - {1\over b} = \sum_{k=a}^{b-1} \left({1\over k} - {1\over k+1} \right) = \sum_{k=a}^{b-1} {1\over k(k+1)} \le \sum_{k=a}^{b-1} {1\over k^2},
%$$
%hence we have
%$$
%{1\over 2^{p+1}} = {1\over 2^p} - {1\over 2^{p+1}} \le \sum_{k=2^p}^{2^{p+1}-1} {1\over k^2}
%$$
%for any $p\ge 0$.

%Since $\Omega$ is non decreasing, we obtain
%$$
%{1\over 2^{p+1}} \log\Omega(2^p) \le \sum_{k=2^p}^{2^{p+1}-1} {1\over k^2} \log\Omega(k), 
%$$
%hence
%$$
%\sum_{\nu\ge 0}{1\over 2^{q+\nu+1}} \log\Omega(2^{q+\nu+2}) \le \sum_{\nu\ge 0} \sum_{k=2^{q+\nu+1}}^{2^{q+\nu+2}-1} {1\over k^2} \log\Omega(k) = \sum_{k\ge 2^{q+1}}{1\over k^2} \log\Omega(k),
%$$
%and we are done. \qed

%\sm{\bf da sistemare e scrivere meglio}

%\sm We would also like to mention here the articles of Zehnder [Z1--3] where one can find the modified Newton method used by R\"ussmann.

%\sm Notice that, when there are no resonances, the function $\omega_f(m)$ defined in \rf{De1.3.7} satisfies
%$$
%|\lambda^Q -\lambda_j|\ge {\omega_f(|Q|)}
%$$
%for each multi-index $Q\in\N$ with $|Q|\ge 2$. 

%Let us then define:

\sm We have the following relation between the R\"ussmann and the reduced Brjuno condition.

\sm\thm{LeRussBrjRed}{Lemma} {\sl Let~$n\ge2$ and let~$\lambda=(\lambda_1, \dots, \lambda_n)\in(\C^*)^n$. If $\lambda$ satisfies R\"ussmann condition, then it also satisfies the reduced Brjuno condition.}

\sm\proof The function $\widetilde\omega_{\lambda_1,\dots\lambda_n}(m)$ defined in Definition \rf{De1.0} satisfies 
$$
\widetilde\omega_{\lambda_1,\dots\lambda_n}(m)^{-1}\le \widetilde\omega_{\lambda_1,\dots\lambda_n}(m+1)^{-1}
$$ 
for all $m\in\N$, and
$$
|\lambda^Q -\lambda_j|\ge \widetilde\omega_{\lambda_1,\dots\lambda_n}(|Q|)
$$
for each $j=1,\dots, n$ and each multi-index $Q\in\N$ with $|Q|\ge 2$ not giving a resonance relative to $j$. Furthermore, by its definition, it is clear that any other function $\Omega\colon\N\to\R$ such that $k\le \Omega(k)\le \Omega(k+1)$ for all $k\in\N$, and satisfying, for any $j=1,\dots n$,
$$
|\lambda^Q -\lambda_j|\ge {1\over \Omega(|Q|)}
$$
for each multi-index $Q\in\N$ with $|Q|\ge 2$ not giving a resonance relative to $j$, is such that 
$$
{1\over\widetilde\omega_{\lambda_1,\dots\lambda_n}(m)}\le \Omega(m)
$$ 
for all $m\in\N$. Hence
$$
\sum_{\nu\ge 0} {1\over p_\nu} \log{1\over\widetilde\omega_{\lambda_1,\ldots,\lambda_n}(p_{\nu+1})}< \sum_{\nu\ge 0} {1\over p_\nu} \log\Omega(p_{\nu+1})
$$
for any strictly increasing sequence of integers~$\{p_\nu\}_{\nu_\ge 0}$ with~$p_0=1$.
Since $\lambda$ satisfies R\"ussmann condition, thanks to Lemma \rf{LeRussmann}, there exists a function $\Omega$ as above such that
$$
\sum_{\nu\ge 0} {1\over s_\nu} \log \Omega(s_{\nu+1}) <+\io,
$$ 
with $\{s_\nu\}$ be defined by $s_\nu := 2^{q+\nu}$, with $q\in\N$, and we are done. \qed

\sm We do not know whether the R\"ussmann condition is equivalent to the reduced Brjuno condition in the multi-dimensional case. As we said, R\"ussmann is able to prove that this is true in dimension one, but to do so he strongly uses the one-dimensional characterization of these conditions via continued fraction.

%%%%%%%%%%%%%%%%%%%%%%%%%%%%%%%%%%%%%%%%%%%%%%%%%%%%%%%%%%%%%%%%%%%%%%%%%%%%%%%%
\vbox{\vskip.75truecm\advance\hsize by 1mm
	\hbox{\centerline{\sezfont References}}
	\vskip.25truecm}\nobreak
%%%%%%%%%%%%%%%%%%%%%%%%%%%%%%%%%%%%%%%%%%%%%%%%%%%%%%%%%%%%%%%%%%%%%%%%%%%%%%%%
\parindent=40pt

\bib{A1} {\sc Abate, M.:} {\sl Discrete local holomorphic dynamics,} in ``Proceedings of 13th Seminar of Analysis and its Applications, Isfahan, 2003'', Eds. S. Azam et al., University of Isfahan, Iran, 2005, pp. 1--32.

\bib{A2} {\sc Abate, M.:} {\sl Discrete holomorphic local dynamical systems,} to appear in ``Holomorphic Dynamical Systems'', G. Gentili, J. Guenot, G. Patrizio eds., Lectures notes in Math., Springer Verlag, 
Berlin, 2009, arXiv:0903.3289v1.

\bib{Ar} {\sc Arnold, V. I.:} ``Geometrical methods in the theory of ordinary differen\-tial equations'', Springer-Verlag, Berlin, 1988.

\bib{Bra} {\sc Bracci, F.:} {\sl Local dynamics of holomorphic diffeomorphisms,} Boll. UMI (8), 7--B (2004), pp. 609--636.

\bib{Brj} {\sc Brjuno, A. D.:} {\sl Analytic form of differential equations,} Trans. Moscow Math. Soc., {\bf 25} (1971), pp. 131--288; {\bf 26} (1972), pp. 199--239.

\bib{D} {\sc Dulac, H.:} {\sl Recherches sur les points singuliers des \'equationes diff\'erentielles}, J. \'Ecole polytechnique II s\'erie cahier IX, (1904), pp. 1--125.

\bib{P} {\sc Poincar\'e, H.:} ``\OE uvres, Tome I'', Gauthier-Villars, Paris, 1928, pp. XXXVI--CXXIX.

\bib{R1} {\sc Raissy, J.}: {\sl Linearization of holomorphic germs with quasi-Brjuno fixed points}, Math. Z., {http://www.springerlink.com/content/3853667627008057/fulltext.pdf}, Online First, (2009).

\bib{R2} {\sc Raissy, J.}: {\sl Simultaneous linearization of holomorphic germs in presence of resonances}, Conform. Geom. Dyn. {\bf 13} (2009), pp 217--224. 

\bib{R3} {\sc Raissy, J.}: {\bf Geometrical methods in the normalization of germs of biholomorphisms}, Ph.D Thesis, Universit\`a di Pisa (2009).

\bib{Re} {\sc Reich, L.}: {\sl Das Typenproblem bei formal-biholomorphien Abbildungen mit anziehendem Fixpunkt}, Math. Ann., {\bf 179} (1969), pp 227--250.

\bib{R\"u1} {\sc R\"ussmann, H.}: {\sl On the convergence of power series transformations of analytic mappings near a fixed point into a normal form}, Preprint I.H.E.S., {\bf M/77/178}, 1977.

\bib{R\"u2}{\sc R\"ussmann, H.}: {\sl Stability of elliptic fixed points of analytic area-preserving mappings under the Brjuno condition}, Ergodic Theory Dynam. Systems, {\bf 22}, (2002), pp 1551--1573.

\bib{Y} {\sc Yoccoz, J.-C.:} {\sl Th\'eor\`eme de Siegel, nombres de Bruno et polyn\^omes quadratiques,} Ast\'erisque {\bf 231} (1995), pp. 3--88.

\bye